\documentclass[12pt]{article}

\usepackage{amsmath}
\usepackage{amssymb}
\usepackage{latexsym}
\usepackage{epsf}
\usepackage{supertabular}
\usepackage{rotating}
\usepackage{multicol}
\usepackage{multirow}
\usepackage{epsfig}
\usepackage{fancyhdr}       

\newif\ifger
\gerfalse

\newtheorem{theorem}{Theorem}[section]

\newtheorem{lemma}[theorem]{Lemma}

\newtheorem{remark}[theorem]{Remark}

\newtheorem{proposition}[theorem]{Proposition}

\title{Flag-transitive automorphism groups of $2$-designs with $\lambda\geq (r, \lambda)^2$ are not product type}

\author{Huiling Li$^1$\ ,  Zhilin Zhang$^2$ and Shenglin Zhou$^3$\footnote{Supported by the National Natural Science Foundation of China (Grant No.11271123). Corresponding author: slzhou@scut.edu.cn}\\
{\small\it 1. School of Mathematical Science, Zhejiang University,}\\
{\small\it  Hangzhou, Zhejiang 310027, PR China}\\
 {\small\it 2. School of Mathematical Science, South China Normal University,}\\
 {\small\it Guangzhou, Guangdong 510631, PR China}\\
{\small\it 3. School of Mathematics, South China University of Technology,}\\
 {\small\it Guangzhou, Guangdong 510640, PR China}}

\begin{document}
\maketitle

\begin{abstract} \quad
In this note we show that a flag-transitive automorphism group $G$ of
a non-trivial $2$-$(v,k,\lambda)$ design with $\lambda\geq (r, \lambda)^2$ is not of product action type. In conclusion, $G$ is a primitive group of affine or almost simple type.

\medskip\noindent
{\bf 2000 Mathematics Subject Classification:} 05B05, 05C25, 20B25

\medskip\noindent
{\bf Key words:}  2-design; flag-transitive; product action

\end{abstract}

\parskip 1ex
\section{Introduction}

A $2$-$(v, k, \lambda)$ design $\cal D$ is a pair $(\Omega, {\cal B})$ with a set $\Omega$ of $v $ points and a set $\cal B$ of blocks such
that each block is a $k$-subset of $\Omega$ and each two distinct points are contained in $\lambda$ blocks.
The design $\cal D$ is nontrivial if $2 < k < v$. Let $b$ denote the number of blocks, then the number of blocks containing each point of $\cal D$ is a constant, denoted by $r$. The positive integers $v, b, r, k, \lambda$ are known as the parameters of $\cal D$.
 An automorphism of a design ${\cal D}=(\Omega,\cal B)$ is a
permutation $\pi$ of $\Omega$ such that $B\in \cal B$ implies
$B^\pi\in\cal B$. The set of all automorphisms of a design forms a group, the (full) automorphism group, $Aut({\cal D})$.
 For a group $G\leq  Aut({\cal D})$, $G$ is said
to be point-primitive if $G$ acts primitively on $\Omega$, and said to be point-imprimitive otherwise.
A flag of $\cal D$ is a point-block pair $(\alpha, B)$ where $\alpha$ is a point and $B$ is a block incident with $\alpha$.
A subgroup $G$ of $Aut({\cal D})$ is said to be flag-transitive if $G$ acts transitively on the set of flags of $\cal D$.

This paper continues to study  the $2$-$(v,k,\lambda)$  designs with  $\lambda\geq (r,\lambda)^{2}$. In this case, any flag-transitive automorphism group $G$
 must act point-primitively on $\cal D$ according to a result of Dembowski \cite[(2.3.7)]{Demb1968}.  In  \cite{Zhou}, Zhou and Zhan studied this type of
 designs which admit a flag-transitive automorphism group $G$.
 They proved the following result.

 \begin{proposition}\label{th1}{\rm \cite[Theorem 1.2]{Zhou}}
Let $\mathcal{D}$ be a $2$-$(v,k,\lambda)$ design with
$\lambda\geq(r,\lambda)^2$. If $G$ is a flag-transitive
automorphism group of $\mathcal{D}$, then $G$ is of affine, almost simple type, or product type with ${\rm Soc}(G)\cong T\times T$, where $T$ is a nonabelian simple group and $G$ has rank $3$.
\end{proposition}

In this paper we focus on the last case of Proposition \ref {th1} in which $G$ has the product action. The main result is:
\begin{theorem} \label{MainTheo1}
Let $\mathcal{D}$ be a $2$-$(v,k,\lambda)$ design with $\lambda\geq(r,\lambda)^2$ which admits a flag-transitive automorphism group $G$. Then $G$ is not of product type.
\end{theorem}

Combing Theorem \ref{MainTheo1} with Proposition \ref{th1}, we obtain the following reduction theorem.

\begin{theorem} \label{cor}
Let $\mathcal{D}$ be a $2$-$(v,k,\lambda)$ design with $\lambda\geq(r,\lambda)^2$ which admits a flag-transitive automorphism group $G$. Then $G$ is of affine or almost simple type.
\end{theorem}

We first give a well known result on 2-designs.

\begin{lemma}{\rm(\cite[Remark, 6.10]{Handbook})}\label{arith}\
Let $\cal D$ be a $2$-$(v, k, \lambda)$ design with parameters $v, b, r, k, \lambda$, then the following hold:
\begin{enumerate}
\item[\rm(i)]\, $r(k-1)=\lambda(v-1)$.
\item[\rm(ii)]\,$vr=bk$.
\item[\rm(iii)]\,$b\geq v$ and $k \leq r$.
\end{enumerate}
\end{lemma}

We finish this section by an elementary lemma. The proof is left to the reader.
\begin{lemma} \label{e} Let $e$ be a positive integer.
\begin{enumerate}
\item[\rm(i)]\, If $e\geq 5$ then $18e^2<2^{2e}$.
\item[\rm(ii)]\,If $e\geq 2$ then $18e^2<3^{2e}$.
\item[\rm(iii)]\, For any  prime $p\geq 5$, $18e^2<p^{2e}$.
\item[\rm(iv)]\,If $e\geq 7$ then $2e^2<2^e$.

\end{enumerate}
\end{lemma}

\section{Proof of Theorem \ref{MainTheo1}}

Assume by the way of contradiction that there exists a $2$-$(v, k, \lambda)$ design $\cal D$ which  satisfies the assumptions of Theorem \ref{MainTheo1}: $G\leq Aut({\cal D})$ has product action on $\Omega=\Delta\times \Delta$, the point set of $\cal D$. We shall denote the points of $\Omega$ by lower case Latin letters, and the points of $\Delta$  by lower case Greek letters. So the point of $\Omega$ has the form $p=(\alpha,\beta)$. We may write
$\Omega=\Delta_1\times \Delta_2$, where $\Delta_1=\Delta_2=\Delta$, and assume that $|\Delta|=\omega\geq 5$  and $\omega$ is odd (see \cite[Section 3.3]{Zhou}) in the following.

By Proposition \ref{th1}, we assume that $T^2\trianglelefteq G\leq T_0\wr \mathbb{Z}_2$, where $T_0$ is 2-transitive on $\Delta$ and $T=Soc(T_0)$ is simple, nonabelian and $v=\omega^2$.
We will prove Theorem \ref{MainTheo1} in three subsections.

\subsection{Some arithmetic results}

\begin{lemma} \label{1}
\begin{enumerate}
\item[\rm(i)] $r^*=2(\omega-1)$ where $r^*=\frac{r}{(r,\lambda)}$.
\item[\rm(ii)] $k=\frac{\lambda^*(\omega+1)}{2}+1$ where $\lambda^*=\frac{\lambda}{(r,\lambda)}$.
\item[\rm(iii)] $b<4v$.
\end{enumerate}
\end{lemma}
{\bf Proof.}\, (i) It has been proved in \cite[Lemma 3.5]{Zhou}.
(ii) By $r^*(k-1)=\lambda^*(v-1)$ and $r^*=2(\omega-1)$, we have $k=\frac{\lambda^*(\omega+1)}{2}+1.$
 (iii) Since $\lambda^{*}\geq (r,\lambda)$, then $$b=\frac{vr}{k}=\frac{2v(\omega-1)(r,\lambda)}{\frac{\lambda^*(\omega+1)+2}{2}}
 <4v\cdot\frac{\omega-1}{\omega+1}\cdot\frac{
(r,\lambda)}{\lambda^{*}}<4v.$$  $\hfill\square$

\begin{lemma} \label{2}
 $k=\displaystyle\frac{(c+d-2)(\omega+1)}{2}+1$ for some positive integers $c$ and $d$.
\end{lemma}
{\bf Proof.}\, It is well known that the suborbits of $G_p$ on $\Omega$ correspond to the orbits of $G$ on $\Omega^{2}=\Omega\times \Omega$. Namely, these orbits are $\Omega_0^{2}=\{(p,p) \,| \, p\in \Omega\},$
$\Omega_1^{2}=\{(p,q) \,| \, p\in \Omega, q\in \Omega_1(p)\}$, $\Omega_2^{2}=\{(p,q) \,| \, p\in \Omega, q\in \Omega_2(p)\},$ where $\Omega_1(p)$ and  $\Omega_2(p)$ are the nontrivial orbits of $G_p$.
Clearly, if $p=(\alpha,\beta)$ and $q=(\gamma,\delta)$, then $(p,q)\in \Omega_1^{2}$ if and only if either $\alpha=\gamma, \beta\neq \delta$, or $\alpha\neq \gamma, \beta=\delta$, and  $(p,q)\in \Omega_2^{2}$ if and only if  $\alpha\neq\gamma$ and $ \beta\neq \delta$.

 Count the pairs $((p,q),B)\in \Omega_1^{2}\times {\mathcal B}$ with $p, q\in B$. As $ |\Omega_1^{2}|=2\omega^2(\omega-1)$ and as each such pair is incident with $\lambda$ blocks we get  $2\omega^2(\omega-1)\lambda$ pairs. If $B$ contains $x$ pairs from $\Omega_1^{2}$ the equation
$2\omega^2(\omega-1)\lambda=bx$ follows.
Similarly count the pairs $((p,q),B)\in \Omega_2^{2}\times {\mathcal B}$ with $p, q\in B$ in two ways, we have the equation
$\omega^2(\omega-1)^2\lambda=by$, where $y$ is the number of pairs from $\Omega_2^{2}$ that contained in $B$. Hence
$y=\frac{x(\omega-1)}{2}$.

Let $B$ be a block of $\cal D$, and $p=(\alpha,\beta)\in B$. Assume that there are $c$ points in $B$ with $\alpha$ as the first component, i.e., $p=(\alpha,\beta),  (\alpha,\beta_1),...,(\alpha,\beta_{c-1})$, and there are $d$ points in $B$  with $\beta$ as the second component, i.e., $p=(\alpha,\beta),  (\alpha_1,\beta),...,(\alpha_{d-1},\beta)$. Then $c$ and $d$ are independent of the choice of $p$ since $G$ is flag-transitive. It is easy to know that there are $c+d-2$ points $q$ in $B$ such that the pair $(p,q)$ is contained in $\Omega_1^{2}$, and  there are $y=\frac{x(\omega-1)}{2}=\frac{(c+d-2)(\omega-1)}{2}$ points $q$ in $B$ such that the pair $(p,q)$ is contained in $\Omega_2^{2}$. From
$(c+d-2)+\frac{(c+d-2)(\omega-1)}{2}=\frac{(c+d-2)(\omega+1)}{2},$
we know that this sum exactly is the number of points of $B$ except $p$, that is $k-1$. The lemma is proved.  $\hfill\square$

\begin{lemma} \label{3}
  $\lambda^*=c+d-2$.
\end{lemma}

{\bf Proof.}
 This follows from Lemmas \ref{1} and \ref{2}. $\hfill\square$

\begin{lemma} \label{4}
 $d=c+1$.
\end{lemma}
{\bf Proof.}\, Without loss of generality, we assume that $c\leq d$. Let
$$B^1=\{\alpha\in \Delta_1\,|\, (\alpha,\beta)\in B\},\ \ \ B^2=\{\beta\in \Delta_2 \,|\, (\alpha,\beta)\in B\}.$$

 Clearly, $|B^1|=\frac{k}{c}$ and $|B^2|=\frac{k}{d}$.
If $c=d$, then $\lambda^*=2(c-1)$ is even, which contradicts the fact that $(r^*,\lambda^*)=1$. If $d\geq c+2$, then by Lemma \ref{2},
$$\frac{k}{c}=\frac{(c+d-2)(\omega+1)}{2c}+\frac{1}{c}>\frac{2c}{2c}(\omega+1)=\omega+1>|\Delta_1|.$$
This contradicts the fact that $|B^1|\leq |\Delta_1|$.
 $\hfill\square$

Note that $c$ and $d=c+1$ are divisors of $k$, we let $k=c(c+1)m$ with $m$ a positive integer. Then by Lemma \ref{2}
\begin{equation} c(c+1)m-\frac{(2c-1)(\omega+1)}{2}=1.
\end{equation}

In the following, we discuss some properties of Equation (1).

\begin{lemma} \label{a} Let $c$ be a positive integer. Then Equation (1) involving the variables $m,\omega$ has positive integer solutions if and only if $c\not\equiv 2 {\pmod 3}$.
\end{lemma}

{\bf Proof.}\, If $c\not\equiv 2 {\pmod 3}$, then $(c(c+1),2c-1)=1$, and there exist integers $x$ and $y$ such that $c(c+1)x+(2c-1)y=1$. For any integer $u$, let $x'=x+u(2c-1)$ and $y'=y-uc(c+1)$. Clearly, $c(c+1)x'+(2c-1)y'=1$. So we can choose $u$ such that $x'$ is a positive integer, and writing $x'=m$ and $y'=-\frac{\omega+1}{2}$, then we have  $c(c+1)m-\frac{(2c-1)(\omega+1)}{2}=1$.
For the converse, suppose that Equation (1) has a positive integer solution $(m,\omega)$. Then $c(c+1)m+(2c-1)\frac{(-\omega-1)}{2}=1$, and hence $(c(c+1),2c-1)=1$ which implies that $c\equiv 0, 1 {\pmod 3}$.  $\hfill\square$

\begin{lemma} \label{b} Let $c\not\equiv 2 {\pmod 3}$  and $(m_0,\omega_0)$ be a positive solution of Equation (1) with $m_0$ minimal. Then the general solution $(m,\omega)$ of Equation (1) has the form:
\begin{enumerate}
\item[(i)] if $c=3t$, then $m=m_0+s\Delta m$, $\omega=\omega_0+s\Delta \omega$, where $m_0=2t+1, \Delta m=6t-1$, $\omega_0=6t^2+6t+1$, $\Delta \omega=18t^2+6t$;
\item[(ii)] if $c=3t+1$, then $m=m_0+s\Delta m$, $\omega=\omega_0+s\Delta \omega$, where $m_0=4t+2, \Delta m=6t+1$, $\omega_0=12t^2+16t+5$, $\Delta \omega=18t^2+18t+4$, where $s$ and $t$ are non-negative integers.
\end{enumerate}
\end{lemma}
{\bf Proof.} \, We only prove part (i), for the proof of part (ii) is similar. Let $c=3t$.
It is easy to check that $(m_0,\omega_0)$ is  a solution of Equation (1) and that $(\Delta m,\Delta \omega)$ is a solution of the equation $$2c(c+1)x-(2c-1)y=0.$$
Thus $(m_0+s\Delta m,\omega_0+s\Delta\omega)$ is a solution of (1) for any integer $s$.
Conversely, if $(m,\omega)$ is a solution of (1), then $(m-m_0,\omega-\omega_0)$ is a solution of the equation of $2c(c+1)x-(2c-1)y=0$. Since $(2c(c+1),2c-1)=1$, $2c(c+1)\mid y$ and $(2c-1)\mid x$, from these we know that every integer solution of (1) has the form $(m_0+s\Delta m, \omega_0+s\Delta \omega)$. Since $6t-1>2t+1$, $m_0=2t+1$ is minimal. $\hfill\square$



\begin{lemma} \label{c} If $(m,\omega)$ is a positive solution of Equation (1), then $2m^2>\omega$.

\end{lemma}

{\bf Proof.}  If $c=3t$, then $m=(2t+1)+s(6t-1)$, $\omega=(6t^2+6t+1)+s(18t^2+6t)$. It is easy to know that $2m^2>\omega$. The proof of the case of $c=3t+1$ is similar.  $\hfill\square$

\begin{lemma} \label{d} If $(m,\omega)$ is a positive solution of Equation (1) for some integer $c\geq 1$, then $(c+1)m<\omega$. Furthermore,

(i)\, if $c\neq 1$, then $(c+1)m<\omega<(c+2)m$;

 (ii)\, if $c=1$, then $\omega=4m-3$ where $m\geq 2$.
\end{lemma}

{\bf Proof.} From Equation (1), we have
$$(c+1)m=\frac{2c-1}{2c}(\omega +1)+\frac{1}{c}
=(\omega+1)-\frac{\omega-1}{2c}.$$
Then $(c+1)m\geq \omega$ if and only if $\frac{\omega-1}{2c}=0$ or 1. Clearly,  $\frac{\omega-1}{2c}\neq 0$ for $\omega\geq 5$. Therefore, $\frac{\omega-1}{2c}=1$ which implies $\omega=2c+1$, and so  $c(c+1)m=(2c-1)(c+1)+1$. It follows that $(c+1)\mid 1$ which is impossible. Thus $(c+1)m< \omega$.

On the other hand, from $(c+1)m=\frac{2c-1}{2c}(\omega +1)+\frac{1}{c}>\frac{2c-1}{2c}\omega$ we know that $\frac{2c(c+1)m}{2c-1}>\omega$. If $c>1$, from $\frac{2c(c+1)}{2c-1}=c+1+\frac{c+1}{2c-1}\leq c+2$ we have $(c+2)m\geq \frac{2c(c+1)m}{2c-1}>\omega$, and so
$(c+1)m<\omega<(c+2)m$.  If $c=1$, by the equation (1), we have $\omega=4m-3$.
$\hfill\square$

As the application of above lemmas, we give the following general result.

\begin{proposition}\label{p2}  Let the notation $\Omega, \Delta$ and $\omega$ be as above, and let $G=(S_{\Delta}\times S_{\Delta})\rtimes \mathbb{Z}_2$, where $S_{\Delta}$ denotes the symmetric group on $\Delta$. If $(m,\omega)$ is a positive integer solution of Equation (1) for some positive integer $c$, then there exists a $2$-$(\omega^2,k,\lambda)$ design $\mathcal D$ and $G$ acts flag-transitively on it.
\end{proposition}

{\bf Proof.}\, Let $G= (S_{\Delta}\times S_{\Delta})\rtimes \mathbb{Z}_2$, then the elements of the subgroup $S_{\Delta}\times S_{\Delta}$ mapping $\Delta_i$ into $\Delta_i$, and the elements which are not in $S_{\Delta}\times S_{\Delta}$ interchange $\Delta_1$ and $\Delta_2$. We use the same symbols such as $\alpha,\beta,...$ to denote the elements of $\Delta_1$ and $\Delta_2$, and $p=(\alpha,\beta)$ denotes the point of $\Omega$ where $\alpha$ is  the first component which means it come from $\Delta_1$, $\beta$ is the second component which means it come from $\Delta_2$. Then $G$ is a permutation group acting on $\Omega$ with rank 3, and $G$
has three orbits on $\Omega^{2}=\Omega\times \Omega$, denoted by $\Omega_0^{2}$, $\Omega_1^{2}$ and $\Omega_2^{2}$ as in Lemma \ref{2}.
Correspondingly, if $p=(\alpha,\alpha)\in \Omega$ is given, then $G_p$, the stabilizer of $p$ in $G$, has 3 orbits
$\Omega_0(p)=\{p\}, \Omega_1(p)=\{(\alpha,\beta), (\beta,\alpha)\,|\, \beta\in \Delta, \alpha\neq\beta\}$ and $\Omega_2(p)=\{(\beta,\gamma)\,|\, \beta,\gamma\in \Delta\setminus \{\alpha\}\}$ of lengths $1$, $2(\omega-1)$ and $(\omega-1)^2$.

Suppose that $c, m,\omega$ satisfy Equation (1) where $\omega=|\Delta|$. By Lemma \ref{d}, we have $(c+1)m<\omega$, it follows that there are  $(c+1)m$ elements in $\Delta_1$:
$$\alpha_1^{(1)}, \alpha_2^{(1)},..., \alpha_{c+1}^{(1)},\alpha_1^{(2)},\alpha_2^{(2)},..., \alpha_{c+1}^{(2)},\alpha_1^{(m)},\alpha_2^{(m)},..., \alpha_{c+1}^{(m)}.$$

Similarly, $cm<(c+1)m<\omega$, there are  $cm$ elements in $\Delta_2$:
$$\beta_1^{(1)}, \beta_2^{(1)},..., \beta_{c}^{(1)},..., \beta_{1}^{(m)},\beta_2^{(m)},..., \beta_{c}^{(m)}.$$

Using these distinct elements, $\alpha_i^{(s)}$ and $\beta_{j}^{(s)}$ where $1\leq i\leq c+1$, $ 1\leq j\leq c$ and $ 1\leq s\leq m$, we construct the points of $\Omega$. Let $p_{ij}^{(s)}=(\alpha_i^{(s)},\beta_{j}^{(s)})$. Then we obtain $c(c+1)m$ points. Put $k:=c(c+1)m$, and
$$B=\{p_{ij}^{(s)}|\,1\leq i\leq c+1, 1\leq j\leq c, 1\leq s\leq m \}.$$
Clearly, $|B|=k$. Let ${\cal B}=\{B^g \,|\, g\in G\}$. Then we assert that the incidence structure ${\cal D}=(\Omega, \cal B)$ is a $2$-$(v,k,\lambda)$ design where $v=\omega^2=|\Omega|$.

In fact, we only need to prove that there exists a positive integer $\lambda$ such that each pair of points $p,q$ of $\Omega$ occurs together in exactly $\lambda$ blocks.

Assume that $p\neq q\in \Omega$, and $p, q$ contained in $\mu$ blocks, say, $B_1, B_2,..., B_{\mu}$. If there exists $g\in G$ such that $p'=p^g, q'=q^g$, then $p', q'$ also contained in $\mu$ blocks, $B_1^g, B_2^g,...,B_{\mu}^g$. Hence, in this case, the number of blocks containing $p$ and $q$ is depending on whether or not $(p,q)$ lie in $\Omega_1^{2}$ or $\Omega_2^{2}$.

If $(p,q)\in \Omega_1^{2}$, $p,q\in B$, where $p=p_{ij}^{(s)}, q=p_{\ell m}^{(s)}$. By the definition of $\Omega_1^{2}$, we have $i=\ell$ and  $j\neq m$, or $j=m$ and $i\neq \ell$. Therefore, when $p$ is fixed, there are $2c-1$ possibilities for $q$ in $B$ such that $(p,q)\in \Omega_1^{2}$, and so among the ordered pairs $(p,q)$ with $p\neq q\in B$ there are $k(2c-1)$ ordered pairs $(p,q)$ belong to $\Omega_1^{2}$. Thus the number of the other pairs is $k(k-1)-k(2c-1)=\frac{k(\omega-1)(2c-1)}{2}.$

Assume that $(p,q)\in \Omega_1^{2}$, and it contained exactly in $\lambda_1$ blocks. Now count the pairs $((p,q),B)\in \Omega_1^{2}\times {\mathcal B}$ with $p, q\in B$ in two ways to get $2v(\omega-1)\lambda_1=bk(2c-1)$, and so $\lambda_1=\frac{bk(2c-1)}{2v(\omega-1)}$.
Similarly, assume that $(p,q)\in \Omega_2^{2}$ which contained exactly in $\lambda_2$ blocks. Count the pairs $((p,q),B)\in \Omega_2^{2}\times {\mathcal B}$ with $p, q\in B$ in two ways, we have the equation  $v(\omega-1)^2\lambda_2=b\frac{k(\omega-1)(2c-1)}{2}$.
Therefore, $$\lambda_2=\frac{bk(2c-1)(\omega-1)}{2v(\omega-1)^2}=\frac{bk(2c-1)}{2v(\omega-1)}=\lambda_1.$$
Putting $\lambda_1=\lambda_2=\lambda$, thus two distinct points $p$ and $q$  are contained in exactly $\lambda$ blocks,  the incidence structure $\cal D$ is a $2$-$(\omega^2,k,\lambda)$ design.

Clearly, $G$ is block-transitive by the structure of $\cal D$. Let
\begin{align*}
 \sigma&=(\alpha_1^{(1)}, \alpha_2^{(1)},...,\alpha_{c+1}^{(1)}),\\
  \tau&=(\beta_1^{(1)}, \beta_2^{(1)},...,\beta_{c}^{(1)}),\\
 \phi&=(\alpha_1^{(1)}, ...,\alpha_{1}^{(m)})(\alpha_2^{(1)}, ...,\alpha_{2}^{(m)})...(\alpha_{c+1}^{(1)}, ...,\alpha_{c+1}^{(m)})(\beta_1^{(1)},...,\beta_1^{(m)})...
 (\beta_c^{(1)},...,\beta_c^{(m)}),
 \end{align*}
then $\sigma\in S_{\Delta_1}$, $\tau\in S_{\Delta_2}$ and $\phi\in S_{\Delta_1}\times S_{\Delta_2}$.
It is easy to know that $\sigma,\tau$ and $\phi$ leave the block $B$ invariant, and so
$\langle \sigma,\tau, \phi\rangle\leq G_B$.
Furthermore, we have $(\alpha_1^{(1)})^{\sigma^{i-1}}=\alpha_i^{(1)}$,  $(\beta_1^{(1)})^{\tau^{j-1}}=\beta_j^{(1)}$ where $1\leq i\leq c+1$ and $1\leq j\leq c$, and  $(p_{11}^{(1)})^{\sigma^{i-1}\tau^{j-1}}=p_{ij}^{(1)}$,
$(p_{ij}^{(1)})^{\phi^{s-1}}=p_{ij}^{(s)}$ where $1\leq s\leq m$, so we conclude that $\langle \sigma,\tau,\phi\rangle$ is transitive on $B$. Thus, $G$ acts flag-transitively on $\cal D$. $\hfill\square$


\begin{remark} In Proposition \ref{p2}, let $c=1$, then $\omega=4m-3$ by Lemma \ref{d}(ii). Also from Lemma \ref{2}, $k=2m$. Thus we have a $2$-$((4m-3)^2,2m,\lambda)$ design admitting a flag-transitive automorphism group with the product action, and the design constructed in Proposition \ref{p2} is a generalization of Davies' example {\rm(}see {\rm\cite[p.53]{Davies})}.
\end{remark}

\subsection{Results on automorphism group}
 Recall that we have proved that there exists a positive integer $c$ such that $k=\frac{(2c-1)(\omega+1)}{2}+1$. For a fixed point $p=(\alpha,\beta)\in B$  with $(p,q)\in \Omega_1^2$ such that $q\in B$, $c$ is the number of points in $B$ with $\alpha$ as  the first component, $d$ is the number of points in $B$ with $\beta$ as  the second component. Note that here $d=c+1$. From this we know that $c$ and  $c+1$ are divisors of $k$, and $k=c(c+1)m$. Hence, $m, \omega$ satisfy Equation (1). By Proposition \ref{p2} we know that  the positive integers $m,\omega$ and $c$ satisfying Equation (1) determine a flag-transitive $2$-$(v,k,\lambda)$ design.


Let  $\cal D$ be a $2$-$(v,k,\lambda)$ design with $\lambda\geq (r,\lambda)^2$. Let
$G\leq Aut({\cal D})$ be flag-transitive which  has the product action on $\Omega=\Delta_1\times \Delta_2$, and $Soc(G)=T^2$ for some nonabelian simple group $T$. So we may write $Soc(G)=T_1\times T_2$ where $T_i\cong T$ for $i=1,2$. Let $A_i=N_{S_{\Delta_i}}(T_i)$. Then $A_i\leq Aut(T_i)$, and $|A_i/T_i|\mid |Out(T_i)|$. It follows that $G\leq (A_1\times A_2)\rtimes \mathbb{Z}_2$. Let
$G_1=(A_1\times A_2)\cap G$, then $|G:G_1|=2$, and $G_1$, as a subgroup of $A_1\times A_2$, is a subdirect product of some subgroups of $A_1$ and $A_2$. Let $P_1=A_1\cap G$, $P_2=A_2\cap G$. Then $P_1\times P_2\unlhd G_1$.
Let
$$
\begin{array}{lll}
Q_1&=&\{x \,|\, x\in A_1, {\mbox{there\,\,exists}}\,\, y\in A_2\,\, {\mbox{such\,\, that}}\, \, (x,y)\in G_1\},\\
Q_2&=&\{y \,|\, y\in A_2, {\mbox{there\,\,exists}}\,\, x\in A_1 \,\, {\mbox{such\,\, that}}\, \, (x,y)\in G_1\}.
\end{array}
$$
Then $Q_1/P_1\cong Q_2/P_2\cong G_1/P_1\times P_2$. Since $T_i\leq P_i\lhd Q_i\leq A_i$  for $i=1,2$, then $Q_i/P_i$ is isomorphic to a section of the group  $A/T$ where $A=N_{S_{\Delta}}(T)\cong A_i$, i.e., there are two subgroups $P, Q$ such that $T\leq P\lhd Q\leq N_{S_{\Delta}}(T)$ with $Q/P\cong Q_i/P_i$. Next we will discuss the structure of $G_B$ and our aim is to obtain the result of ``$m\leq |A/T|$''.

1)\, $G_B\not\leq P_1\times P_2$.

Assume for contrary that $G_B\leq P_1\times P_2$. Then corresponding to the block-stabilizer $G_B$, we define the following 4 subgroups:
$$
\begin{array}{lll}
L_1&=&\{x\,|\, x\in P_1, (x,1)\in G_B\},\\
L_2&=&\{y\,|\, y\in P_2, (1,y)\in G_B\},\\
M_1&=&\{x\,|\, x\in P_1, {\mbox{there\,\,exists}}\,\, y\in P_2\,\, {\mbox{such\,\, that}}\, \, (x,y)\in G_B\},\\
M_2&=&\{y\,| \, y\in P_2, {\mbox{there\,\,exists}}\,\, x\in P_1 \,\, {\mbox{such\,\, that}}\, \, (x,y)\in G_B\}.
\end{array}
$$
Then $L_1\times L_2\leq G_B\leq M_1\times M_2$, and $G_B/L_1\times L_2\cong M_1/L_1\cong M_2/L_2$. Clearly, $L_i=G_B\cap P_i$, and $M_i$ is the projection of $G_B$ onto $P_i$, and hence we have the following subgroups series:
$$ G>G_1\geq P_1\times P_2\geq M_1\times M_2\geq G_B\geq L_1\times L_2,$$
where $|G:G_1|=2$. Now we show that $|P_1\times P_2: M_1\times M_2|\geq \omega^2=v$.

In fact, since  $|P_1\times P_2:M_1\times M_2|=|P_1:M_1||P_2:M_2|$, we only prove that $|P_1:M_1|\geq \omega$, the other case is similar. Recall that $B^1$ is defined in Lemma \ref{4} and it is the orbit of $M_1$ with length $\frac{k}{c}$. If $B^1=\Delta_1$, then $\omega=\frac{k}{c}=\frac{(2c-1)(\omega+1)}{2c}+\frac{1}{c}=(\omega+1)-\frac{\omega-1}{2c}$ and so $\omega=2c+1$. It follows that $d=(c+1)\nmid k$, a contradiction. So that $B^1\subset\Delta_1$.
  Now suppose that $|B^1|\leq \omega-1$. From $(B^1)^{M_1}=B^1$, we have $M_1\leq (P_1)_{B^1}$. Let ${\mathcal D}_1 $ be the incidence structure with the point set $\Delta_1$ and the block set $\{(B^1)^g\,|\, g\in P_1\}$. Then ${\cal D}_1$ is a 2-design, and $|P_1:M_1|\geq |P_1: (P_1)_{B^1}|=b_1\geq \omega$ by Lemma \ref{arith}(iii), where $b_1$ is the number of blocks of ${\mathcal D}_1$.

If $M_1\times M_2\neq G_B$, then by above subgroups series we have $b=|G:G_B|\geq 2\cdot v\cdot 2=4v$. This contradicts Lemma \ref{1}(iii). Thus $G_B=M_1\times M_2$, and then
 $M_i=L_i$ for $i=1,2$.

Both $B^1$ and $B^2$  may be regarded as subsets of $\Delta$ and $P_1\cong P_2:\equiv P$ a permutation group on $\Delta$, we see that the incidence structure $\mathcal{D}_1$ with $\Delta$ as the point set and  $\{(B^1)^g\,|\,g\in P\}$ as the block set is a $2$-$(\omega, \frac{k}{c},\lambda_1)$ design and the incidence structure $\mathcal{D}_2$ with $\Delta$ as the point set  and $\{(B^2)^h\,|\,h\in P\}$  as the block set is a $2$-$(\omega, \frac{k}{c+1},\lambda_2)$ design. Let $b_i, r_i$ be the parameters of $\mathcal{D}_i$ for $i=1,2$.

So the block $B$ is just a subset of $\Omega$ with  the form $\{(\alpha, \beta)\,|\, \alpha\in B^1, \beta \in B^2\}$. In the following  we use $(B^1, B^2)$ to denote this subset.

Let $B=(B^1, B^2)$. Since $P_1\times P_2$ is a subgroup of the automorphism group $G$ of $\cal D$, for any $g\in P_1$ and $ h\in P_2$, $((B^1)^g,(B^2)^h)$ is also a block of $\cal D$. On the other hand, there exists an element of $G$ which interchanges two components of points $(\alpha,\beta)$ of $\Omega$, the subsets $((B^2)^h,(B^1)^g)$, where $g\in P_1$ and $ h\in P_2$, are also blocks of $\cal D$. It is easy to see the above two kinds of  blocks constitute the block set of $\cal D$. Let $\alpha,\beta,\gamma,\delta$ be four elements of $\Delta$, then $p=(\alpha,\beta), q=(\alpha,\delta)$ and  $s=(\gamma,\delta)$ are three distinct points of $\Omega=\Delta_1\times \Delta_2$. Let $(X,Y)$ be any block of $\cal D$. Therefore $p$ and $q\in (X,Y)$ if and only if $\alpha\in X$ and $\beta,\delta\in Y$. There are $r_1$ possibilities for $X$ and $\lambda_2$ possibilities for $Y$, so that $p, q$ are contained in $r_1\lambda_2+r_2\lambda_1$ blocks of $\cal D$ because $p,q$ are contained in $r_1\lambda_2$ blocks of the form $((B^1)^g,(B^2)^h)$ and are contained in $r_2\lambda_1$ blocks of the form $((B^2)^h,(B^1)^g)$. On the other hand, $p$ and $s\in  (X,Y)$ if and only if $\alpha,\gamma\in X$ and $\beta,\delta\in Y$. So that $p, s$ are contained in $2\lambda_1\lambda_2$ blocks of $\cal D$. However, $r_1\lambda_2+r_2\lambda_1>2\lambda_1\lambda_2$ for $r_i>\lambda_i(i=1,2)$, and hence $\cal D$ is not a 2-design. Thus, $G_B\not\leq P_1\times P_2$ as asserted.

2)\, The structure of $G_B$.

Consider the product $G_B(P_1\times P_2)=\langle G_B, P_1\times P_2\rangle$. Since $G_B\not\leq P_1\times P_2$, then $P_1\times P_2\leq G_B(P_1\times P_2)\leq G_1$, and
$G_B(P_1\times P_2)/P_1\times P_2 \cong G_B/G_B\cap(P_1\times P_2)$. It follows that
$|G_B(P_1\times P_2):G_B|=|P_1\times P_2:G_B\cap(P_1\times P_2)|$. Similarly, we construct 4 subgroups corresponding to $G_B\cap(P_1\times P_2)$ as follows:
$$
\begin{array}{lll}
L_1'&=&\{x\,|\, x\in P_1, (x,1)\in G_B\cap(P_1\times P_2)\},\\
L_2'&=&\{y\,|\, y\in P_2, (1,y)\in G_B\cap(P_1\times P_2)\},\\
M_1'&=&\{x\,|\, x\in P_1, {\mbox{there\,\,exists}}\,\, y\in P_2\,\, {\mbox{such\, that}}\, \, (x,y)\in G_B\cap(P_1\times P_2)\},\\
M_2'&=&\{y\,|\, y\in P_2, {\mbox{there\,\,exists}}\,\, x\in P_1 \,\, {\mbox{such\, that}}\, \, (x,y)\in G_B\cap(P_1\times P_2)\}.
\end{array}
$$
It is easy to know that $L_i'=L_i$ for $i=1,2$. We have proved in 1) that $G_B\not\leq P_1\times P_2$, and so $M_i\not\leq P_i$. But, by the definition of $M_i'$, we always have $M_i'\leq P_i$. The group $G_B\cap (P_1\times P_2)$ satisfies
$$L_1\times L_2 \leq G_B\cap(P_1\times P_2)\leq M_1'\times M_2',$$ and $$G_B\cap(P_1\times P_2)/L_1\times L_2\cong M_1'/L_1\cong M_2'/L_2.$$
Consider the subgroups series
$$G>G_1\geq G_B(P_1\times P_2)\geq P_1\times P_2\geq M_1'\times M_2'\geq G_B\cap(P_1\times P_2). $$
 If  $M_1'\times M_2'\neq G_B\cap(P_1\times P_2)$, similarly we have $|P_1\times P_2:M_1'\times M_2'|\geq \omega^2=v$. Now by the facts $|G_B(P_1\times P_2):G_B|\geq 2v$ and $|G:G_1|=2$, we get $b=|G:G_B|\geq 4v$, a contradiction. Hence, $G_B\cap(P_1\times P_2)=M_1'\times M_2'=L_1\times L_2$.

3)\, Since $G_B\cap(P_1\times P_2)=L_1\times L_2$, we have $G_B/L_1\times L_2\cong G_B/G_B\cap(P_1\times P_2)\cong G_B(P_1\times P_2)/P_1\times P_2\leq G_1/P_1\times P_2$,  which is isomorphic to a section of $A/T$, and hence $$|G_B:L_1\times L_2|\mid |A/T|.$$

4)\, The structure of $B$.

Recall that if $(\alpha,\beta)\in B$, then $B^1$ is the  orbit of $M_1$ containing $\alpha$,
$B^2$ is the  orbit of $M_2$ containing $\beta$. Since $L_i\lhd M_i$, the group $L_i$ acts $1/2$-transitively on $B^i$, so $B^1$ is the join of orbits of $L_1$, say, $\varphi_1$, $\varphi_2$,..., $\varphi_s$ which has the same length $c'$,  and $B^2$ is the join of orbits, say $\psi_1, \psi_2,..., \psi_t$ of $L_2$ which has the same length $d'$. Let $(\alpha,\beta)\in B$, and $\alpha\in \varphi_i$, $\beta\in\psi_j$. Then for any $\alpha'\in \varphi_i$, $\beta'\in \psi_j$, the points $(\alpha', \beta')$ are all in $B$, and we called the set $\varphi_i\times \psi_j$ is the $(\varphi,\psi)$-set containing $(\alpha, \beta)$. By above assumptions we know that each $(\varphi,\psi)$-set contains $c'd'$ points, and $B$ is formed by finitely many $(\varphi,\psi)$-sets. Therefore, $k=c'd'm'$ for some positive integer $m'$.

We have also defined two integers $c$ and $d$ and proved that $d=c+1$ in Subsection 2.1: Let $p=(\alpha,\beta)\in B$, there are $c$ points in $B$  with $\alpha$ as the first component, there are $d$ points in $B$  with $\beta$ as the second component. By the flag-transitivity of $G$, $c$ and $d$ are independent of the choice of $p=(\alpha,\beta)$.   Since there are $c'$ points in $\varphi_i\times \psi_j$ have the same second component, $d'$ points have the same first component, then $c'\leq c$, $d'\leq d=c+1$, and hence $m\leq m'$. Since  $G_B$ acts transitively on the $(\varphi,\psi)$-sets which formed $B$, and the kernel is $L_1\times L_2$, then
 $m'\mid |G_B/L_1\times L_2|$. Combining this with $m\leq m'$ and $|G_B/L_1\times L_2|\mid |A/T|$ implies
 \begin{equation}
 m\leq |A/T|.
 \end{equation}

\subsection{The final contradiction}

Since in our situation the nonabelian simple groups $T$ which can occur as minimal normal subgroups of 2-transitive groups $T_0$ of degree $\omega$ is also $2$-transitive on $\Delta$ in all cases except $T=PSL(2,8)$ of degree 28 (1-transitive)(\cite[Notes 2]{Cameron81}),  we list $T$, $\omega$ and $|Out(T)|$ in the following (refer to \cite{Cameron81,Kantor85,LPS90}). Note that in the list, $\omega$ is odd, $q=p^e$ where $p$ is a prime, and $A/T=Out(T)$ except Case (2).
\begin{enumerate}
\item[(1)] $T= A_\omega$, $\omega \geq 5$, $|Out(T)|=2$;
\item[(2)] $T = PSL(d, q)$, $d \geq 2$, $q$ is a prime power and $\omega = \frac{q^d-1}{q-1}$, where $(d, q)\neq$ $(2,2), (2,3)$, $|Out\, T|=2(d,q-1)e$ if $d\geq3$,  or $(d,q-1)e$ if $d=2$;
\item[(3)] $T = PSU(3, q)$, $q \geq 3$ is a prime power and $\omega = q^3 + 1$, $|Out\, T|=2(3,q+1)e$;
\item[(4)] $T = Sz(q)$, $\omega = q^2+1$, $|Out(T)|=e$;
\item[(5)] $T=PSL(2,11)$, $\omega = 11$, $|Out(T)|=2$;
\item[(6)] $T=A_7$, $\omega=15$, $|Out(T)|=2$;
\item[(7)]$T= M_{\omega}$, $\omega = 11, 23$, $|Out(T)|=1$.
\end{enumerate}

 In what follows, we analyse each of these possible cases separately. 

{\bf Case (1):} Here $A/T=Out(T)$.  By $m\leq |A/T|$, $m=1$ or $2$. Then from Lemma \ref{c}, $2m^2>\omega$ and $\omega\geq 5$, we get  $\omega=5$ or $7$ and $m=2$. Now Lemma \ref{d}, $2(c+1)<\omega\leq 2(c+2)$ implies  $(\omega,c)=(5,1)$ or $(7, 2)$. By Lemma \ref{a}, $c\neq 2$. If $(\omega,c)=(5,1)$,
then $d=2$, and Lemma \ref{3} implies $\lambda^*=2c-1=1$. It follows that $\lambda=(r,\lambda)$.
But, the condition $\lambda\geq (r,\lambda)^2$ gives $\lambda=(r,\lambda)=1$. This can be ruled  by \cite[Theorem]{Zieschang}.

{\bf Case (2):} Let $T = PSL(d, q)$ with $\omega = \frac{q^d-1}{q-1}$, where $(d, q)\neq$ $(2,2), (2,3)$, $|Out(T)|=2(d,q-1)e$ if $d\geq3$,  or $(d,q-1)e$ if $d=2$.

Note that, in this case, $T$ is 2-transitive on $\Delta$, and so $T$ is the projective special linear group $PSL(d,q)$ of $d$-dimensional vector space $V$ over the field $\mathbb{F}_q$, and $\Delta$ is the set of all 1-dimensional subspaces of $V$.  Let $\sigma\in Aut(T)$, then $\sigma \in N_{S_{\Delta}}(T)$ if and only if $\sigma$ transforms any 1-dimensional subspace of $V$ to a 1-dimensional subspace. However, since the graph automorphism $\gamma$ of $T$ transforms the stabiliser of any 1-dimensional subspace to  the stabiliser of a $(d-1)$-dimensional subspace, then $\gamma\not\in A=N_{S_{\Delta}}(T)$ if $d\geq 3$. Hence in this case $|A/T|=(d,q-1)e$.

We first show that $d\leq 4$. For the contrary, if $d\geq5$, then $\omega=\frac{q^d-1}{q-1}>q^{d-1}\geq q^4$. Also by the fact $m\leq |A/T|$, we have
$m\leq (d,q-1)e\leq (q-1)e<eq$.
Then $2m^2<2e^2q^2$.  It is easily known by induction that $2e^2< 2^{2e}$ for any positive integer $e$, and so $2e^2<2^{2e}\leq p^{2e}=q^2$. Therefore, $2m^2<2e^2q^2<q^4<\omega$, contradicting Lemma \ref{c}.

{\bf Subcase (2.1):}  $d=4$.

Here $\omega=q^3+q^2+q+1$. Since $\omega$ is odd, then $q=2^e$, and so $m\leq |A/T| =(4,q-1)e=e$. Then $2m^2\leq 2e^2<q^2<\omega$, contradicting Lemma \ref{c}.

{\bf Subcase (2.2):}  $d=3$.

Here $\omega=q^2+q+1>q^2$ and $|A/T|=(3,q-1)e\leq 3e$. It follows from the inequality (2) and Lemma \ref{c} that $p^{2e}<2m^2\leq 18e^2$.

By Lemma \ref{e}(i)-(iii), except cases $(p,e)\in \{(2,1),(2,2), (2,3), (2,4),(3,1)\}$, we have $2m^2\leq 18e^2< p^{2e}=q^2<\omega$, which contradicts to Lemma \ref{c}.

Now we deal with the remaining five pairs $(p,e)$ in the following.
For these pairs $(p,e)$, we list $T$, $|A/T|$, $\omega$ and $m$ in Table \ref{tab:1}.
\begin{table}[!h]
\centering
\caption{$d=3$}
\label{tab:1}
\begin{tabular}{ccccc}
\hline
{\sc Line}&$T$&$|A/T|$&$\omega$&$m$ \\
\hline
1&$PSL(3,2)$&1&7&{\mbox{non-existence}}\\
2&$PSL(3,4)$&6&21&$4,5,6$\\
3&$PSL(3,8)$&3&73& {\mbox{non-existence}}\\
4&$PSL(3,16)$&12&273&$12$\\
5&$PSL(3,3)$&1&13& {\mbox{non-existence}}\\
\hline
\end{tabular}
\end{table}
Note that Column 5 gives the values of $m$ by the following method. Since $2m^2>\omega$, then $m>\sqrt {\frac{\omega}{2}}$. On the other hand, we have $m\leq |A/T|$, and so $[\sqrt{\frac{\omega}{2}}]+1\leq m\leq |A/T|$. Then using the information of Columns 3 and 4 we obtain the possible values of $m$. Note that {\sc Lines} 1, 3 and 5 cannot occur because there is no positive integers $m$ satisfying above inequality.

{\sc Line} 2: Here $\omega=21$, $m=4,5$ or 6.
If $m=4$ or 5, then Equation (1) changed to $2c^2-9c+5=0$ or $5c^2-17c+10=0$.  But each equation has no positive integer solutions,  a contradiction.
If $m=6$,  then Equation (1) changed to $3c^2-8c+5=0$,  and so $c=1$.
Now Lemmas  \ref{3} and \ref{4} imply $\lambda^*=2c-1=1$. So that $\lambda=(r,\lambda)$. The condition $\lambda\geq (r,\lambda)^2$ implies $(r,\lambda)=\lambda=1$. Therefore, by \cite[Theorem]{Zieschang}, this case cannot occur. (In fact, the case $(c,m)=(1,6)$ also can be ruled out by Lemma \ref{b}(ii).)

{\sc Line} 4: Here $\omega=273$. Then by Equation (1), we have $12c(c+1)=137(2c-1)+1$, i.e., $6c^2-131c+68=0$. This equation has no positive integer solutions, a contradiction.
 Thus, $d\neq 3$.

{\bf Subcase (2.3):}  $d=2$.

Here  $\omega=q+1>q$, and $q=2^e$ since $\omega$ is odd. Then we have $m\leq |A/T|=e$ and $2m^2<2e^2$. By Lemma \ref{e} (iv), if $e\geq 7$ then $2e^2<2^e$, and combing it with $2m^2<2e^2$, we get $2m^2<q<\omega$, a contradiction. Hence $e=2,3,4,5$ or 6.

Similarly, we have Table \ref{tab:2}.
\begin{table}[!h]
\centering
\caption{$d=2$}
\label{tab:2}
\begin{tabular}{ccccc}
\hline
{\sc Line}&$T$&$|A/T|$&$\omega$&$m$ \\
\hline
1&$PSL(2,4)$&2&5&2\\
2&$PSL(2,8)$&3&9&3\\
3&$PSL(2,16)$&4&17&3, 4\\
4&$PSL(2,32)$&5&33&5\\
5&$PSL(2,64)$&6&65&6\\
\hline
\end{tabular}
\end{table}

 Now for each line in Table \ref{tab:2}, we solve the equation (1) on variable $c$ and obtained the desired contradiction.

 {\sc {Line}} 1: Here $2c(c+1)=3(2c-1)+1$, and so $c^2-2c+1=0$. Hence $c=1$, and $\lambda^*=2c-1=1$. So that $\lambda=(r,\lambda)$. The condition $\lambda\geq (r,\lambda)^2$ implies $(r,\lambda)=\lambda=1$. This can be ruled by \cite[Theorem]{Zieschang}.

{\sc Line} 2: Here $3c(c+1)=5(2c-1)+1$, and so $3c^2-7c+4=0$. Hence $c=1$. As above we have $(r,\lambda)=\lambda=1$. This can be ruled by \cite[Theorem]{Zieschang}.


{\sc Lines} 3-5: Here, for each value of $m$, we get a quadratic equation. But it is easy to know that each equation has no positive integer solutions, a contradiction. Therefore, $d\neq 2$.

{\bf Case (3):} Here $A/T=Out\,T$, and so $|A/T|=|Out\, T|=2(3,q+1)e\leq 6e$. Then from $m\leq |A/T|$ we have $m\leq 6e$, and $2m^2<72e^2$. On the other hand,
by the inequality $2m^2>\omega$ we get $2(6e)^2>\omega >q^3$. Since $\omega$ is odd, let $q=2^e\geq 3$. Then
$9e^2>2^{3(e-1)}$, which gives $e=2$ or 3.

If $e=2$, then $\omega=65$. Then by Lemma \ref{c}, we have $2m^2> 65$. It follows that $m\geq 6$, which contradicts the fact $m \leq 2(3, q+1)e=4$.

If $e=3$,  then $\omega=513$, and from $m\leq |A/T|$ we have $m\leq 2(3, q+1)e=18$. On the other hand, by Lemma \ref{c}, we have $2m^2> 513$, and so $m\geq 17$.  Hence $m=17$ or 18. However, the equation (1) changed to $17c^2-497c+256=0$ or $9c^2-248c+128=0$ respectively,  which has no positive integer solutions, this is a contradiction.

{\bf Case (4):} Here $A/T=Out\,T$, then by $m\leq |A/T|$, we have $m\leq e$. It follows from $2m^2>\omega$ that $2e^2>2^{2e}$, and this inequality does not holds for any $e\geq 1$, which is a contradiction.

{\bf Cases (5)-(6):} Here $\omega = 11$ or 15, and $|A/T|=2$. Then by $m\leq |A/T|$, we have $m=1$ or $2$, which contradicting Lemma \ref{c}.

{\bf Case (7):} Here $\omega = 11$ or 23, and $|A/T|=1$. Then $m=1$, contradicting Lemma \ref{c}.

This complete the proof of Theorem \ref{MainTheo1}. $\hfill\square$

\end{document}